\documentclass[12pt]{amsart}

\usepackage{amssymb}

\usepackage{graphicx}

\newtheorem{thm}{THEOREM}[section]
\newtheorem{prop}[thm]{PROPOSITION}
\newtheorem{lemma}[thm]{LEMMA}
\newtheorem{cor}[thm]{COROLLARY}
\newtheorem{defn}[thm]{\em Definition}

\theoremstyle{definition}

\theoremstyle{remark}
\newtheorem{remark}[thm]{Remark}

\numberwithin{equation}{section}

\def\C{\mathbb{C}}
\def\R{\mathbb{R}}

\def\CP{\mathbb{CP}}

\def\N{\mathbb{N}}

\def\proof{\noindent{\em{Proof}.\ }}
\begin{document}

\title[$1$-convex non-embeddable threefolds]{Some examples of $1$-convex non-embeddable threefolds}

\author{Giovanni Bassanelli}
\address{(G. Bassanelli) Dipartimento di Matematica,
Universit\`a degli Studi di Parma, Parco Area delle Scienze 53/A,
I-43100 Parma, Italy} \email{giovanni.bassanelli@unipr.it}

\author{Marco Leoni}
\address{(M. Leoni) Coll\`ege Louis Nuc\'era,
171 Route de Turin, 06300 Nice, France} \email{marco.leoni2003@libero.it}

\thanks{Research partially supported by the GNSAGA of the INdAM.}


\dedicatory{}


\begin{abstract} We construct a family of  $1$-convex threefolds, with exceptional curve $C$ of type
$(0,-2)$, which are not embeddable in $\C^m \times \CP_n$. In order to show that they are not
K\" ahler we exhibit a real $3$-dimensional chain $A$ whose boundary is the complex curve $C$.
 \end{abstract}

\maketitle

\section{Introduction}

In general a $1$-convex manifold with $1$-dimensional exceptional set is embeddable, that is it
can be realized as an embedded subvariety of $\C^m \times \CP_n$, for suitable $m$ and $n$.
The only  possible exceptions are given by the following theorem:

\begin{thm}\label{I} {\em (see Theorem 3 in \cite{C})} Let $X$ be a non-embeddable, 
$1$-convex manifold whose exceptional set $C$ has dimension $1$. Then $dim_{\C}X=3$
and  $C$ has an irreducible component which is a rational curve of type $(-1,-1)$, $(0,-2)$
or $(1,-3)$.
\end{thm}

\noindent As regards the existence of the quoted exceptions,
in \cite{C} there is an example of type $(-1,-1)$.
 In \cite{V} (p. 242 B) there is an example of type $(0,-2)$, but the argument 
is dubious (see \cite{C}, Remark 4, but see also \cite{V2}). For the case $(1,-3)$ nothing is known. 

In fact the first two cases are easier:  as it is well known (see \cite{L}), there is a family $\{W_k\}_{k\in \N^*}$  
of fiber bundles over $\CP_1$ which give a local model: this means that if $X$ is a $1$-convex threefold whose
exceptional set is a smooth rational curve $C$ of type $(-1,-1)$ (respectively $(0,-2)$) then there is a 
neighbourhood of  $C$ biholomorphic to a neighbourhood of the null section of $W_1$ (resp. of 
$W_k$, for a suitable $k \ge 2$) (see Definition \ref{A1} and Proposition \ref{A}). Moreover the 
 sequence of the normal bundles associated to $C$  is
$$\underbrace{ (0,-2), \dots, (0,-2) }_{k-1},(-1,-1).$$

Starting from this remark we shall show:

\begin{thm}\label{II} For every integer $k \ge 1$ there is a non-embeddable $1$-convex threefold 
$\tilde X_k$  whose exceptional set is a smooth rational curve $C$ for whose the sequence of
normal bundles is
$$\underbrace{ (0,-2), \dots, (0,-2) }_{k-1},(-1,-1).$$
In particular $N_{C|\tilde X_1}=\mathcal O(-1)\oplus \mathcal O(-1)$ and, for $k \ge 2$,
$N_{C|\tilde X_k}=\mathcal O(0)\oplus \mathcal O(-2)$.
\end{thm}

For $k=1$ we get the above quoted Coltoiu's example. 
We point out that our construction is explicit and elementary; in order to see that $\tilde X_k$ is not
K\" ahler we shall exhibit a real $3$-chain $A$ whose boundary is the exceptional curve $C$.

\bigskip
\section{The proof of Theorem \ref{II}}

\begin{defn}\label{A1} {\em Let $k \ge 1$ be an integer. The equations :
\begin{equation}\label{a}\left\{ \begin{array}{l}
w=\frac{1}{z}\\
y_1=z^2x_1+zx_2^k\\y_2=x_2
\end{array} \right.\end{equation}
in the coordinates $(z,x_1,x_2)$ and $(w,y_1,y_2)$, define a fiber bundle on $\CP_1$, 
with fiber $\C^2$, which will be denoted by $W_k$.}
\end{defn}

As we just said in the Introduction, these manifolds $W_k$ are local models for $1$-convex threefolds,
as the following Proposition states.

\begin{prop}\label{A} (see \cite{P}, p. 234) In any  $1$-convex threefold whose exceptional set is a smooth rational curve
$C$ of type $(-1,-1)$ (resp. $(0,-2)$) there exists a neighbourhood of $C$ biholomorphic to a neighbourhood
of the null section of $W_1$ (resp. of $W_k$, for a suitable $k \ge 2$).
\end{prop}

There is a  geometrical description of these threefolds:

\begin{prop}\label{B} Let $k \ge 1$ be an integer. Let $N_k \stackrel{f_k}{\longrightarrow}\C^4$
be the blow-up with center at the complex smooth surface 
\begin{equation}\label{a0}S_k:=\{z\in \C^4;z_1-iz_2=z_3-z_4^k=0\}.\end{equation}
 Then $W_k$ is the strict transform of the hypersurface
\begin{equation}\label{a1}Y_k:=\{z\in \C^4;z_1^2+z_2^2+z_3^2-z_4^{2k}=0\}\end{equation}
 whose only one singular point is the origin
$P_k=0$; the null section of $W_k$ is $C_k=f_k^{-1}(P_k)$. 
\end{prop}
\proof It is enough to follow the outline of \cite{P}, Example 2.14. \qed

\bigskip

Now we shall investigate in more detail this geometric construction and build a suitable
commutative diagram
\begin{equation}\label{d}\begin{array}{ccccccccc}N_0&\stackrel{h_1}{\longrightarrow}&N_1&\stackrel{h_2}{\longrightarrow}&
\dots&\stackrel{h_{k-1}}{\longrightarrow}&N_{k-1}&\stackrel{h_k}{\longrightarrow}&N_k\\
\downarrow\scriptstyle {f_0}& &\downarrow\scriptstyle {f_1}&&&&\downarrow\scriptstyle {f_{k-1}}
&&\downarrow\scriptstyle {f_k}\\
M_0&\stackrel{g_1}{\longrightarrow}&M_1&\stackrel{g_2}{\longrightarrow}&
\dots&\stackrel{g_{k-1}}{\longrightarrow}&M_{k-1}&\stackrel{g_k}{\longrightarrow}&M_k=\C^4
\end{array}\end{equation}

\medskip

\noindent{\bf Step 1.} {\em Applying the desingularization 
process to the hypersurface $Y_k \subset \C^4:=M_k$ we the sequence}
$$M_0\stackrel{g_1}{\longrightarrow} M_1 
\stackrel{g_2}{\longrightarrow}
\dots \stackrel{g_{k-1}}{\longrightarrow} M_{k-1} \stackrel{g_k}{\longrightarrow} M_k=\C^4$$
More precisely this sequence  is defined by induction:  $M_{k-1} \stackrel{g_k}{\longrightarrow} M_k=\C^4$ 
is the blow-up with center 
$P_k:=0$. Then  define the chart \linebreak
$(U_{k-1};u_1,\dots,u_4)$ of $M_{k-1}$ saying that in 
these coordinates
the map $g_k$ has the following equations:
\begin{equation}\label{b}
\left\{ \begin{array}{ll}
z_j=u_ju_4, & \textrm{for}\  j=1,2,3\\
z_4=u_4 &
\end{array} \right.
\end{equation}
Denoting by $Y_{k-1}\subset M_{k-1}$ the strict transform of $Y_k$, we get that the only  singular point of 
$Y_{k-1}$ is $P_{k-1}:=0\in U_{k-1}$ and
\begin{equation}\label{c}Y_{k-1}\cap U_{k-1}=\{u_1^2+u_2^2+u_3^2-u_4^{2(k-1)}=0\}.\end{equation}
 Comparing \eqref{c} with \eqref{a1}, we see that we can iterate the process: the map 
$M_{j-1} \stackrel{g_j}{\longrightarrow}M_j$ is the blow-up of center $P_j$ and $Y_{j-1}$ is the strict transform 
of $Y_j$. Finally, since
$$Y_0\cap U_0=\{u_1^2+u_2^2+u_3^2-1=0\}$$
$Y_0$ is smooth, so that the process ends.

\medskip
We need the following lemma:

\begin{lemma}\label{C} Let $S$ be a smooth complex surface in a complex $4$-fold M and let $P \in S$. 
There is the following commutative diagram:
$$\begin{array}{ccc}N'&\stackrel{h}{\longrightarrow}&N\\
\downarrow\scriptstyle {f'}& &\downarrow\scriptstyle {f}\\
M'&\stackrel{g}{\longrightarrow}&M
\end{array}$$
where: $g$ is the blow-up with center $P$, $S'$ is the strict transform of $S$ in $M'$, $f$ (resp. $f'$) is
 the blow-up of center $S$ (resp. $S'$),  $h$ is the blow-up with center the curve $C:=f^{-1}(P)$.
\end{lemma}

\proof The problem is local near $P$, thus we may assume that $M=\C^4$ and that $S$ is a plane. 
Choosing $S=\C^2 \times\{0\}$
the direct computation is easier. \qed

\bigskip
Now, recalling Proposition \ref{B}, we finish our construction: 
\medskip

\noindent{\bf Step 2.}  {\em Define $S_{j-1}$ as the strict transform of $S_j$ by
means of the map $M_{j-1} \stackrel{g_j}{\longrightarrow}M_j$, $1 \le j \le k$. 
Let 
$N_j \stackrel{f_j}{\longrightarrow} M_j$ be the blow-up of center $S_j$, $0 \le j \le k-1$.
Moreover let $C_j:=f_j^{-1}(P_j)$ and  $N_{j-1} \stackrel{h_j}{\longrightarrow} N_j$ be the blow-up
with center $C_j$, $1 \le j \le k$.

Then the diagram \eqref{d} is commutative.}

\proof   By means of  Lemma \ref{C} it is enough to check that $P_j \in S_j$, $j=0, \dots, k$.
But, as noted above, $P_j=0 \in U_j$ and using the chart $U_j$, it is straighforward to check that
\begin{equation}\label{e}S_j\cap U_j=\{u_1-iu_2=u_3-u_4^j\}.\qed\end{equation}

\begin{cor}\label{E} Let $X_j$ be the strict transform of $Y_j$ by means of the map 
$N_j \stackrel{f_j}{\longrightarrow}M_j$, $0 \le j \le k$. Considering  restrictions
of  maps, we get, from \eqref{d} the following commutative diagram:
\begin{equation}\label{f}\begin{array}{ccccccccc}X_0&\stackrel{h_1}{\longrightarrow}&X_1&\stackrel{h_2}{\longrightarrow}&
\dots&\stackrel{h_{k-1}}{\longrightarrow}&X_{k-1}&\stackrel{h_k}{\longrightarrow}&X_k=W_k\\
\downarrow\scriptstyle {f_0}& &\downarrow\scriptstyle {f_1}&&&&\downarrow\scriptstyle {f_{k-1}}
&&\downarrow\scriptstyle {f_k}\\
Y_0&\stackrel{g_1}{\longrightarrow}&Y_1&\stackrel{g_2}{\longrightarrow}&
\dots&\stackrel{g_{k-1}}{\longrightarrow}&Y_{k-1}&\stackrel{g_k}{\longrightarrow}&Y_k
\end{array}\end{equation}
Moreover:

\begin{itemize}
\item[\em (i)] $X_j$ is smooth, the rational curve $C_j$ (which is the center of center of $h_j$) is contained in $X_j$ and
there is a neighbourhood of $C_j$ in $X_j$ biholomorphic to a neighbourhood
of the null section of $W_j$, $1 \le j \le k$;
\item[\em (ii)] $C_j=h_j(C_{j-1})$, $j \ge 2$;
\item[\em (iii)] $X_0 \stackrel{f_0}{\longrightarrow} Y_0$ is a biholomorphism;
\item[\em (iv)] the exceptional divisor $E_0$ of $h_1$ is biholomorphic to $\CP_1 \times \CP_1$
and the induced map $E_0 \simeq \CP_1 \times \CP_1  \stackrel{h_1}{\longrightarrow} C_1 \simeq \CP_1$
is one of the two canonical projections.
\end{itemize}
\end{cor}
\proof  The diagram \eqref{f} is well defined, because diagram \eqref{d} is commutative.
Comparing equations \eqref{a1} and \eqref{c}, \eqref{a0} and \eqref{e}
 we may apply Proposition \ref{B} for $j=1, \dots, k$. Thus we get $W_k=X_k$ and,
for $1 \le j \le k-1$, $W_j=X_j \cap f_j^{-1}(U_j)$ and $C_j$ is its null section.
By  Proposition \ref{A}, $N_{C_1|X_1}=(-1,-1)$, while $N_{C_j|X_j}=(0,-2)$,
for $j \ge 2$. Thus the exceptional divisor $E_{j-1}$ of $h_j$ is a rational ruled surface: $E_0=F_0$
(this proves {\em (iv)}) and $E_j=F_2$, for $j \ge 1$. Now the curve $C_j$ is not a fiber of $F_2$, 
otherwise $N_{C_j|X_j}=(0,-1)$, thus $C_j$ is a section of $E_j \simeq F_2$; this shows {\em (ii)}.
Finally, since $S_0$ and $Y_0$ are smooth, $X_0 \stackrel {f_0}{\longrightarrow}Y_0$ is a biholomorphism.
\qed

\bigskip

\begin{remark}\label{E1} In the rational ruled surface $F_2$ there is only one curve $C$ 
with negative self-intersection: $C._{F_2}C=-2$. Since $C_j$ is not a fiber of $E_j$, 
from the exact sequence
$$0 \to \mathcal O (C_j._{_{E_j}} C_j) \to N_{C_j|X_j} \to \mathcal O (C_j._{_{X_j}}E_j) \to 0$$
it follows easily that $C_j$ is the curve of $E_j$ with negative self-intersection; this means that
the sequence $X_0 \to \dots \to X_k$ is the sequence of the blow-ups associated to the curve
$C_k$.
\end{remark}

\bigskip
Let us state the 
following elementary result

\begin{lemma}\label{F} Let $\mathcal Q:=\{z \in \CP_3;z_0^2+z_1^2+z_2^2-z_3^2=0\}$. 
Every line of $\mathcal Q$ has a real point.
\end{lemma}
\proof  Let $r \subset \mathcal Q$ 
be a line and let $P \in r $. If $P$ is not real, consider the line $s$ passing 
through $P$ and $\overline P$. Now $s$ is a real line and $s \cap \R^4$ is external to the real sphere 
$\mathcal Q\cap \R^4$, therefore there are exactly two
 planes passing through $s$  tangent to $\mathcal Q$ and the tangent points are real.
One of these planes must be the plane $\alpha$ defined by the 
lines $r$ and $s$ (infact $\alpha \cap \mathcal Q$ contains $r$ and thus is a degenerate conic), therefore 
is tangent to $\mathcal Q$ in a real point $Q$ which must belong to $r$. \qed 
\bigskip

By means of the detailed description of the map 
$X_k \stackrel {f_k}{\longrightarrow}Y_k$ given in \eqref{f},  the following statement
is a simple corollary.

\begin{cor}\label{G} Let $B:=\{z \in Y_k \cap \R^4;z_4>0\}$ and $A:=f_k^{-1}(B)$. Then $A$ is a
real threefold with boundary $\partial A=C_k$.
\end{cor}

\proof Let $D:=(g_k \circ \dots \circ g_1)^{-1}(B)$. Since the diagram \eqref{f} is commutative, 
$A=h_k \circ \dots \circ h_1(f_0^{-1}(D))$. From \eqref{b} it follows that $g_k^{-1}(B) \subset U_{k-1}$, 
and iterating this argument we get that $D \subset U_0$, more precisely
$$D=\{x\in \R^4;x_1^2+x_2^2+x_3^2-1=0,\ x_4>0\}.$$
Therefore the boundary $\partial D=\{x\in \R^4;x_1^2+x_2^2+x_3^2-1= x_4=0\}$ is contained
in $\{z\in \C^4;z_1^2+z_2^2+z_3^2-1= z_4=0\}=E_0 \cap U_0$. By Corollary \ref{E}(iv) 
the fibers of the map
$E_0 \simeq f_0^{-1}(E_0) \stackrel{h_1}{\longrightarrow}C_1$are given by one of the two family of lines
of the quadric $E_0$. By  Lemma \ref{F} each of these lines intersect $\partial D$, therefore
$h_1 (f_0^{-1}(\partial D))=C_1$. Thus from Corollary \ref{E}(ii) $\partial A=C_k$. \qed

\bigskip
In order to obtain our example $\tilde X_k$
we must perturb $Y_k$ outside the origin.

\begin{lemma}\label{6a} For every fixed integer $k \ge 1$ there exist an integer $N > k$ and
$0 < \varepsilon  \le 1$ such that the origin is the only singular point of the hypersurface
\begin{equation}\label{7}\tilde Y_k:=\{z_1^2+z_2^2+z_3^2-z_4^{2k}+
\varepsilon(z_1^{2N}+z_2^{2N}+z_3^{2N}+z_4^{2N})=0\}.
\end{equation}

\end{lemma}

The equations $w_j:=z_j(1+\varepsilon z_j^{2N})^{1/2}$, $1 \le j \le 3$ and 
$w_4=z_4(1-\varepsilon z_4^{2N})^{1/2k}$ define a biholomorphic map  
$V \stackrel{\Phi}{\longrightarrow} \tilde V$ between two neighbourhood of the origin
of $\C^4$. Thus we can define $\tilde X_k$ gluing $f_k^{-1}(V)\cap X_k$ and 
$\tilde Y_k \setminus \{0\}$ by means of $\Phi$. The maps 
$\tilde X_{j-1} \to \tilde X_j$ and 
$\tilde Y_{j-1} \to \tilde Y_j$ are defined as above since nothing
is changed near $P_j$ and $C_j$, while the  maps $\tilde X_{j} \to \tilde Y_j$ are defined
by a gluing process (these maps are not blow-ups).
\begin{prop}\label{7a} $\tilde X_k$ is not embeddable.
\end{prop}

\proof Let $\tilde B:=\R^4 \cap \tilde Y_k \cap \{x_4 > 0\}$. From \eqref{7} it follows that $\tilde B$
is relatively compact. From the definition of $\Phi$ it follows that $\tilde B \cap \tilde V = \Phi^{-1}(B \cap V)$;
thus the preimage $\tilde A$ of $\tilde B$ is again a $3$-chain with boundary $C_k$. Hence
the exceptional curve $C_k=\partial \tilde A$ is a boundary in $\tilde X_k$, which is not K\"ahler. \qed

\bigskip
\bibliographystyle{amsalpha}

\end{document}